\documentclass[12pt]{amsart}

\newif\ifpdfstoll
\ifx\pdfoutput\undefined
  \pdfstollfalse
\else
  \pdfoutput=1
  \pdfstolltrue
\fi

\ifpdfstoll
  \usepackage[pdftex]{epsfig}
  \usepackage[pdftex,colorlinks=true,%
              pdftitle={Rational points on curves},%
              pdfauthor={Michael Stoll}]{hyperref}
  
\else
  \usepackage{epsfig}
  
\fi

\usepackage{amssymb}  
\usepackage{latexsym} 
\usepackage{comment}

\usepackage[all]{xy}

\DeclareFontEncoding{OT2}{}{} 
\newcommand{\textcyr}[1]{%
 {\fontencoding{OT2}\fontfamily{cmr}\fontseries{m}\fontshape{n}\selectfont #1}}

\newcommand{\Sha}{{\mbox{\textcyr{Sh}}}}

\newcommand{\Z}{{\mathbb Z}}
\newcommand{\Q}{{\mathbb Q}}
\newcommand{\R}{{\mathbb R}}
\newcommand{\F}{{\mathbb F}}

\newcommand{\BP}{{\mathbb P}}

\newcommand{\To}{\longrightarrow}

\newcommand{\Sel}{\operatorname{Sel}}

\newcommand{\rank}{{\text{rank}}}
\newcommand{\tors}{{\text{tors}}}

   {\begin{list}{}{\settowidth{\labelwidth}{#1}%
                   \setlength{\leftmargin}{\labelwidth}%
                   \addtolength{\leftmargin}{\labelsep}}}%
   {\end{list}}

\newtheorem{Theorem}{Theorem}

\newtheorem{Corollary}[Theorem]{Corollary}
\newtheorem{Conjecture}[Theorem]{Conjecture}

\theoremstyle{definition}
\newtheorem{Definition}[Theorem]{Definition}
\newtheorem{Example}[Theorem]{Example}
\newtheorem{Examples}[Theorem]{Examples}

\newtheorem{Remarks}[Theorem]{Remarks}

\newtheorem{Problem}[Theorem]{Problem}

\begin{document}

\title{Rational points on curves}

\author{Michael Stoll}
\address{Mathematisches Institut,
         Universit\"at Bayreuth,
         95440 Bayreuth, Germany.}
\email{Michael.Stoll@uni-bayreuth.de}

\date{16 February, 2010}


\begin{abstract}
  This is an extended version of an invited lecture I gave at the
  Journ\'ees Arithm\'etiques in St.~\'Etienne in July~2009.
  
  We discuss the state of the art regarding the problem of finding
  the set of rational points on a (smooth projective) geometrically
  integral curve~$C$ over~$\Q$. The focus is on practical aspects of
  this problem in the case that the genus of~$C$ is at least~$2$,
  and therefore the set of rational points is finite.
\end{abstract}

\maketitle



\section{Introduction}

\subsection{The Problem}

Let $C$ be a geometrically integral curve defined over~$\Q$. We stick to~$\Q$
here for simplicity. In principle, we can replace $\Q$ by an arbitrary number
field. In practice, however, many of the necessary algorithms are only implemented
for~$\Q$, and even when they are available for more general number fields,
the computations are usually much more involved. We consider the following
problem.

\begin{Problem} \label{TheProblem}
  Determine $C(\Q)$, the set of \emph{rational points} on~$C$.
\end{Problem}

We observe that a curve and its smooth projective model
only differ in a computable finite set of points (coming from points at
infinity and from singularities). Therefore we lose nothing if we assume
that $C$ is \emph{smooth and projective}.


\subsection{The Structure of $C(\Q)$}

Before we consider curves of `higher genus' more specifically, let us
recall what is known about the structure of the set~$C(\Q)$ in general.
There is a trichotomy, depending on the \emph{genus}~$g$ of~$C$, which
is the most important geometric (or even topological, if we think of~$C$
as a Riemann surface) invariant of the curve. This exemplifies the belief
that ``Geometry determines arithmetic'' --- the structure of the set of
rational points on a variety should only depend on its geometry.

We have the following three cases.

\begin{itemize}\addtolength{\itemsep}{2mm}
  \item $g = 0$\,: \\
        In this case, we either have $C(\Q) = \emptyset$ (this is always
        possible), or else if there is a point $P_0 \in C(\Q)$,
        then $C$ is isomorphic over~$\Q$ to the projective line~$\BP^1$.
        Any such isomorphism will give us a parameterization of~$C(\Q)$
        in terms of rational functions in one variable. Probably the best-known
        example is the unit circle $x^2 + y^2 = 1$, whose points can be rationally
        parameterized in the following way.
        \[ t \longmapsto \Bigl(\frac{1-t^2}{1+t^2}, \frac{2t}{1+t^2}\Bigr) \]
        As $t$ runs through $\BP^1(\Q) = \Q \cup \{\infty\}$, its image
        runs through all the rational points on the unit circle. So such
        a parameterization gives us a finite description of the set~$C(\Q)$.
  \item $g = 1$\,: \\
        We either have $C(\Q) = \emptyset$, or else if there is a point
        $P_0 \in C(\Q)$, then $(C, P_0)$ is an \emph{elliptic curve}.
        So $C$ has a geometrically defined structure as an abelian group
        with~$P_0$ as its origin. This implies that $C(\Q)$ is also an
        abelian group with origin~$P_0$. Mordell~\cite{Mordell} has shown that
        $C(\Q)$ is \emph{finitely generated}. (Weil~\cite{Weil} has extended
        this to elliptic curves and, more generally, Jacobian varieties over
        arbitrary number fields.) In particular, we can describe
        $C(\Q)$ by listing generators of this group.
  \item $g \ge 2$\,: \\
        This is the case of `higher genus'. Mordell~\cite{Mordell} has
        conjectured, and Faltings~\cite{Faltings} has proved that the set
        $C(\Q)$ is always \emph{finite}. In particular, we can describe
        $C(\Q)$ by simply listing the finitely many points.
\end{itemize}

We see that in each case, there is a finite description of~$C(\Q)$. The
precise version of Problem~\ref{TheProblem} above therefore asks for
an algorithm that provides this description.

Before we consider the higher genus case in detail, let us give a short
discussion of the other two cases.


\subsection{Genus Zero}

If $C$ is a smooth projective geometrically integral curve of genus~0 (over
any field~$k$), then $C$ is isomorphic to a smooth \emph{conic}. If we can
compute in~$k$, then we can find an explicit such isomorphism. This can
be done by computing a basis of the Riemann-Roch space of an anticanonical
divisor; the map to~$\BP^2$ given by this basis provides the desired isomorphism.

Now let $k = \Q$ again.
Like all quadrics, conics~$C$ satisfy the \emph{Hasse Principle}:
If $C(\Q) = \emptyset$, then $C(\R) = \emptyset$
or $C(\Q_p) = \emptyset$ for some prime~$p$, where $\Q_p$ is the field
of $p$-adic numbers. For future reference, we make the followinf definition.

\begin{Definition}
  The curve $C$ has points \emph{everywhere locally}, if $C(\R) \neq \emptyset$
  and $C(\Q_p) \neq \emptyset$ for all primes~$p$.
\end{Definition}

The Hasse Principle then states that a curve that has points everywhere locally
must also have rational points.

Let us assume that $C$ is given by a ternary quadratic form with integral
coefficients. Then $C(\Q_p) \neq \emptyset$ whenever $p$ does not divide
the discriminant of the quadratic form (we use the fact that smooth conics
over finite fields always have rational points, plus Hensel's Lemma to lift
to a $p$-adic point). So there are only finitely many primes to check in
addition to~$C(\R)$.
(One should note, however, that in general one has to \emph{factor} the
discriminant, which can be difficult.) For each given prime, Hensel's Lemma
gives us an upper bound for the $p$-adic precision needed. So the check
whether $C$ has points everywhere locally
reduces to a finite computation. Therefore we can
\emph{decide} if $C(\Q)$ is empty or not. This is still true for a
number field in place of~$\Q$.

If $C(\Q) \neq \emptyset$ and we know the `bad' primes (those dividing
the discriminant), then there is an efficient procedure that exhibits
a point $P_0 \in C(\Q)$, see for example~\cite{Simon}. This can be
seen as a `minimization' process that finds a $\Q$-isomorphic curve
with good reduction at all primes, followed by a `reduction' process
based on lattice basis reduction~\cite{LLL} that brings the curve into
the standard form $y^2 = xz$, which has some obvious points. This
last (reduction) part of the procedure has, to my knowledge, not yet
been generalized to arbitrary number fields in a satisfactory way.

Given $P_0 \in C(\Q)$, we can easily compute an isomorphism $\phi : C \to \BP^1$
by projecting away from~$P_0$. The inverse of~$\phi$ then provides us with
the desired parameterization of~$C(\Q)$.


\subsection{Genus One}

For curves of positive genus, the Hasse Principle no longer holds in general.
So there is no easy way to check if the curve has rational points or not.
If we cannot find a rational point, but $C$ has points everywhere locally,
then we can try to use a \emph{descent} computation. For $n \ge 2$,
\emph{$n$-descent} consists in computing a finite number of $n$-coverings
of~$C$ such that each of these $n$-coverings has points everywhere locally
and every rational point on~$C$ is the image of a rational point on one
of the $n$-coverings. An \emph{$n$-covering} is a morphism $\pi : D \to C$ 
of curves over~$\Q$ that
over~$\bar{\Q}$ is isomorphic to the multiplication-by-$n$ map $E \to E$,
where $E$ is $C/\bar{\Q}$ considered as an elliptic curve. In principle,
this computation is possible for every~$C$ and every~$n$ over every
number field. In practice however, this is feasible only in a few cases.
\begin{itemize}\addtolength{\itemsep}{1mm}
  \item $y^2 = \text{quartic in $x$}$ and $n = 2$~\cite{Cassels,MSS};
  \item intersections of two quadrics in~$\BP^3$ and $n = 2$~\cite{Stamminger};
  \item plane cubics and $n = 3$ (this is a current PhD project).
\end{itemize}

If the finite set of relevant $n$-coverings turns out to be empty, this
proves that $C(\Q) = \emptyset$. If we assume that Shafarevich-Tate groups
of elliptic curves do not contain nontrivial infinitely divisible elements
(this assumption is weaker than the standard conjecture that $\Sha(E/\Q)$ is
finite),
then it follows that if $C(\Q) = \emptyset$, then there must be an~$n$
such that there are no $n$-coverings of~$C$ with points everywhere locally.
This means that we can, at least in principle, verify that $C$ does not
have rational points.

On the other hand, if $C$ does have rational points, then their preimages
on suitable $n$-coverings tend to be `smaller' and can therefore be found
more easily by a search. So $n$-descent on~$C$ serves two purposes: it allows us
to show that no rational points exist, but it can also help us find a rational
point.

It should be noted that if a curve of genus~1 has infinitely many rational
points, the smallest point can be exponentially large in terms of the coefficients
of the defining equations. (This comes from the corresponding property of
generators of the group of rational points on an elliptic curve.) This
phenomenon is what can make life rather hard when we try to find the rational
points on a curve of genus~1.


\subsection{Elliptic Curves}

We now assume that we have found a rational point $P_0$ on our curve~$C$
of genus~1. Then, as mentioned above, $(C, P_0)$ is an \emph{elliptic curve},
which we will denote~$E$. By Mordell's Theorem we know that $E(\Q)$ is a
finitely generated abelian group; our task is now to find explicit \emph{generators}
of this group. By the structure theorem for finitely generated abelian groups,
we have
\[ E(\Q) = E(\Q)_{\tors} \oplus \Z^r \,, \]
where $E(\Q)_{\tors}$ is the finite subgroup of~$E(\Q)$ consisting of all
elements of finite order. This finite subgroup is easy to find.
The hard part is to determine the \emph{rank}~$r = \rank\, E(\Q)$.

We can use $n$-descent again. When we apply it to an elliptic curve, the
set of $n$-coverings with points everywhere locally has a natural group
structure; this group is the \emph{$n$-Selmer group of~$E$}. Its order
gives an upper bound for the size of $E(\Q)/n E(\Q)$, from which we can
deduce an upper bound for~$r$. As before, this computation is always
possible in principle (see~\cite{StollIHP}).
In practice, $n$-descent on an elliptic curve
over~$\Q$ is currently feasible for $n = 2, 3, 4, 8$; the case
$n = 9$ is current work. See~\cite{CFOSS,CFS} for a detailed description.
In some cases, we can use what is known about the conjecture of Birch
and Swinnerton-Dyer. If the conductor of~$E$ is not too large, we can
compute the values of the $L$-series of~$E$ and its derivatives at~$s=1$
to sufficient precision. If $L(E,1) \neq 0$, then $r = 0$, and if
$L'(E,1) \neq 0$, then $r = 1$~\cite{Kolyvagin}.

A search for independent points in~$E(\Q)$ gives a lower bound on~$r$.
However, generators may be very large. In the same way as for
general curves of genus~1, descent can help us to find them.
When $r = 1$, \emph{Heegner points} can be used if the conductor of~$E$
is sufficiently small.

\begin{Example}
  (See~\cite{CFS}.) The group $E(\Q)$, where
  \[ E : y^2 = x^3 + 7823 \,, \]
  is infinite cyclic and generated by the point
  \begin{align*}
     P = \Bigl(&\frac{2263582143321421502100209233517777}{11981673410095561^2}, \\[3mm]
               &\frac{186398152584623305624837551485596770028144776655756}%
                    {11981673410095561^3}\Bigr) \,.
  \end{align*}
  This point was found by a 4-descent on~$E$. The Heegner point method is
  not feasible here, because the conductor of~$E$ is $2^4 \cdot 3^3 \cdot 7823^2$,
  which is a bit too large.
\end{Example}

A discussion of how one can try to find the set of rational points on
an elliptic curve, or more generally, on a genus~1 curve, would provide
enough material for at least one book. But this is a different story
and will be told at another occasion.


\section{Checking Existence of Rational Points}

We now turn to curves of higher genus, meaning $g \ge 2$. The first question
we would like to answer is whether there are any rational points on the
curve~$C$ or not. 


\subsection{Finding Points}

If $C(\Q)$ is nonempty, we can usually find a rational
point by search. This is because (in contrast to the case of genus~1)
we expect the rational points to be fairly \emph{small}.
The following is a consequence of Vojta's Conjecure; see Su-Ion Ih's
thesis~\cite{Ih}.

\begin{Conjecture}
  If ${\mathcal C} \to B$ is a family of higher-genus curves, then
  there are constants $\gamma$ and~$\kappa$ such that
  \[ H_{\mathcal C}(P) \le \gamma H_B(b)^\kappa \qquad
      \text{for all $P \in {\mathcal C}_b(\Q)$}
  \]
  if the the fiber ${\mathcal C}_b$ is smooth.
\end{Conjecture}
Here $H_B$ denotes a (non-logarithmic) height on the base~$B$, and
$H_{\mathcal C}$ is a suitable height function on~${\mathcal C}$.

If $C$ is hyperelliptic, one can use the {\tt ratpoints} program~\cite{Ratpoints}
for the point search.

\begin{Examples}
  Consider a curve
  \[ C : y^2 = f_6 x^6 + \dots + f_1 x + f_0 \]
  of genus 2, with $f_j \in \Z$.
  Then the conjecture says that there are $\gamma$ and~$\kappa$ such that
  the $x$-coordinate $p/q$ of any point $P \in C(\Q)$ satisfies
  \[ |p|, |q| \le \gamma \max\{|f_0|, |f_1|, \dots, |f_6|\}^\kappa \,. \]
  
  In~\cite{BruinStollExp}, we consider curves of genus~2 as above such that
  $f_j \in \{-3, -2, \dots, 3\}$. We will call these curves
  \emph{small genus 2 curves} in this paper.
  If such a curve has rational points,
  then there is one whose $x$-coordinate is $p/q$ with $|p|,|q| \le 1519$.
  In fact, for all but two such curves (up to isomorphism), we even
  have $|p|,|q| < 80$. On the other hand, the largest point known on one
  of these curves (which is very likely the largest point there is) has
  height~$209\,040$, which indicates that $\gamma$ and/or~$\kappa$ cannot
  be too small.
\end{Examples}

So usually we can assume that we know all the points in~$C(\Q)$. In particular,
if we are unable to find a rational point on~$C$, it is reasonable to suspect
that there are indeed no rational points. The problem now is to \emph{prove}
this fact in some way.


\subsection{Local Points}

One approach that we can try is to check if $C$ has points everywhere
locally. As before, this can be done by a finite computation, which is
efficient modulo the determination of the `bad' primes. This usually comes
down to factoring some kind of discriminant. In addition to the bad primes,
one also has to look at small primes. The reason is that smooth curves
of genus~$g$ may fail to have $\F_p$-points when $p$ is small relative to~$g$.
(By the Weil bounds, we have $\#C(\F_p) \ge p + 1 - 2g \sqrt{p}$, so there
will be $\F_p$-points whenever $p+1 > 2g\sqrt{p}$.)

\begin{Example} (Poonen-Stoll~\cite{PoonenStoll})
  About 84--85\% of all curves of genus~2 have points everywhere locally.
  This percentage is a \emph{density}: we consider all genus~2 curves
  of the form $y^2 = f(x)$ with $f = f_6 x^6 + \dots + f_1 x + f_0 \in \Z[x]$
  such that $\max\{|f_j|\} \le N$, and determine the proportion $\alpha_N$
  of curves with points everywhere locally. Then $\lim_{N\to\infty} \alpha_N$
  exists and has approximately the value given above. Convergence seems
  to be rather fast, compare the data given at the end of~\cite{BruinStoll2Desc}.
\end{Example}

The counterpart to this result is the following conjecture.

\begin{Conjecture} 
  0\% of all curves of genus~2 have rational points.
\end{Conjecture}

In fact, heuristic considerations suggest the following. Let $\beta_N$
be the proportion of curves of size up to~$N$ that possess rational points
(similarly to~$\alpha_N$ above). Then $\beta_N \ll N^{-1/2}$.
See~\cite{StollDensity} for details and some experimental data.

This indicates that checking for points everywhere locally will usually
not suffice to prove that $C(\Q) = \emptyset$: the Hasse Principle is
quite likely to fail.

\begin{Example} (Bruin-Stoll~\cite{BruinStollExp})
  Among the 196\,171 isomorphism classes of small genus 2 curves,
  there are 29\,278 that are counterexamples to the Hasse Principle.
\end{Example}


\subsection{Descent Again}

So we need another method of attack. One possibility is again \emph{descent}.
We find a covering $\pi : D \to C$ (more precisely, an unramified covering of
smooth projective geometrically integral curves that over~$\bar{\Q}$ is
a Galois covering). As before in the genus~1 case, this covering has finitely
many \emph{twists} $\pi_\xi : D_\xi \to C$ such that $D_\xi$ has points
everywhere locally.

\begin{Example}
  Consider a hyperelliptic curve
  \[ C : y^2 = g(x) h(x) \]
  with $\deg g$, $\deg h$ not both odd. Then
  \[ D: \quad u^2 = g(x)\,,\quad v^2 = h(x) \]
  is an unramified $\Z/2\Z$-covering of~$C$ with $\pi : D \to C$ given
  by $(x,u,v) \mapsto (x, uv)$. Its twists are
  \[ D_d : \quad d u^2 = g(x)\,,\quad d v^2 = h(x)\,, \qquad
      d \in \Q^\times/(\Q^\times)^2 \,.
  \]
  Every rational point on~$C$ lifts to one of the twists, since $g(x)$ must
  have some square class~$d$. If $g$ and~$h$ have integral coefficients
  and $p$ is a prime divisor of~$d$ (we assume $d$ to be a squarefree integer),
  then $D_d$ does not have $p$-adic points unless $g$ and~$h$ have a common
  root mod~$p$. This is the case only when $p$ divides the resultant of $g$ and~$h$.
  So we see that only finitely many of the twists~$D_d$ can have
  points everywhere locally.
\end{Example}

The idea of descent goes back to Fermat (`descente infinie'). The statement
that only finitely many twists are relevant is a variant of a result due
to Chevalley and Weil~\cite{ChevalleyWeil}, see Theorem~\ref{Thm:Descent} below.

Here is a concrete example.

\begin{Example}
  Consider the genus~2 curve 
  \[ C : y^2 = -(x^2 + x - 1)(x^4 + x^3 + x^2 + x + 2) =: f(x) \,. \]
  $C$ has points everywhere locally. This can be seen by observing that
  \[ f(0) = 2\,, \quad f(1) = -6\,, \quad f(-2) = -3 \cdot 2^2\,,
     \quad f(18) \in (\Q_2^\times)^2 \quad\text{and}\quad 
     f(4) \in (\Q_3^\times)^2.
  \]
  The first three values show that $C(\R) \neq \emptyset$ and that
  $C(\Q_p) \neq \emptyset$ for all $p \neq 2, 3$; the last two fill the
  remaining gaps. 

  The relevant twists of the obvious $\Z/2\Z$-covering are among
  \[ d\,u^2 = -x^2 - x + 1\,, \qquad d\,v^2 = x^4 + x^3 + x^2 + x + 2 \]
  where $d$ is one of $1$, $-1$, $19$ or $-19$, since the resultant of the two
  factors is 19.
  If $d < 0$, the second equation has no solution in~$\R$;
  if $d = 1$ or $19$, the pair of equations has no solution over~$\F_3$. 
  This is because the first equation implies
  that $x \bmod 3$ is one of $0$ or~$-1$, whereas the second equation implies
  that $x \bmod 3$ is one of $1$ or~$\infty$.

  So there are no twists with points everywhere locally, and therefore
  $C(\Q) = \emptyset$.
\end{Example}

The general result is as follows.

\begin{Theorem}[Descent Theorem] \label{Thm:Descent}
  Let $\pi : D \to C$ be an unramified covering that is geometrically Galois.
  Its twists $\pi_\xi : D_\xi \to C$ are parameterized by $\xi \in H^1(\Q, G)$
  (a Galois cohomology set), where $G$ is the Galois group of the covering.
  We then have the following:
  \begin{itemize}\addtolength{\itemsep}{2mm}
    \item $C(\Q) = \displaystyle\bigcup_{\xi \in H^1(\Q, G)}
            \pi_\xi\left(D_\xi(\Q)\right)$.
    \item $\Sel^{\pi}(C) := \left\{\xi \in H^1(\Q, G) : 
            \text{$D_\xi$ has points everywhere locally}\right\}$ \\
          is finite (and computable).
  \end{itemize}
\end{Theorem}

\begin{Definition}
  In the situation of the Descent Theorem, we call $\Sel^{\pi}(C)$
  the \emph{Selmer set} of~$C$ w.r.t.~$\pi$.
\end{Definition}

\begin{Corollary}
  If we find $\Sel^{\pi}(C) = \emptyset$, then $C(\Q) = \emptyset$ as well.
\end{Corollary}

This follows from the two statements in the theorem, since $D_\xi(\Q)$ is
empty unless $D_\xi$ has points everywhere locally.


\subsection{Abelian Coverings}

In principle, we can use this approach with any covering of~$C$ in the
above sense. However, in practice it is easier to restrict to a special
kind of coverings.

\begin{Definition}
  A covering $\pi : D \to C$ as above is \emph{abelian} if its Galois group~$G$
  is abelian.
\end{Definition}

The reason why abelian coverings are especially useful comes from the
following fact (which is a result of `Geometric Class Field Theory';
see~\cite{SerreGCFT} for details). We let $J$ denote the \emph{Jacobian variety}
of~$C$. We assume for simplicity that there is an embedding $\iota : C \to J$.
(This means that there is a divisor class of degree~1 on~$C$ that is defined
over~$\Q$, i.e., stable under the action of the absolute Galois group
of~$\Q$. If we can show that there is no such divisor class, then it
follows that $C$ does not have rational points, since a rational point
would provide us with a suitanble divisor class.)

Then all abelian coverings of~$C$ are obtained from \emph{$n$-coverings} of~$J$:
\begin{equation} \label{CovDiagram}
   \xymatrix{ D \ar[r] \ar[d]_{\pi}
               & X \ar[d] \ar@{-->}[r]^{\cong/\bar{\Q}}
               & J \ar[dl]^{\cdot n} \\
              C \ar[r]^{\iota} & J
            }
\end{equation}
Here $X \to J$ is an $n$-covering of~$J$, meaning that there is an isomorphism
of $X$ with~$J$ over~$\bar{\Q}$ that makes the triangle in the diagram
commute, and $\pi : D \to C$ is the pull-back of $X \to J$ under~$\iota$.

We call such a covering $D \to C$ an \emph{$n$-covering} of~$C$;
the set of all $n$-coverings with points everywhere locally
is denoted $\Sel^{(n)}(C)$ and called the \emph{$n$-Selmer set} of~$C$.
Every abelian covering of~$C$ can be extended to an $n$-covering for
some~$n$. Therefore the Jacobian gives us a handle on all the
abelian coverings of~$C$. The process of computing the set $\Sel^{(n)}(C)$
is called an \emph{$n$-descent on~$C$}.


\subsection{Computing $n$-Selmer Sets in Practice}

In practice, computing $\Sel^{(n)}(C)$ is usually quite hard, even though
it is possible in principle. The most difficult obstacle is that the
computation requires arithmetic information like ideal class groups
and unit groups for number fields that can be rather large. About the only fairly
general situation where the fields involved are manageable is the computation
of the 2-Selmer group of a \emph{hyperelliptic curve} $C : y^2 = f(x)$.
In this case, the relevant number fields are those generated by a root
of each irreducible factor of~$f$. 
This is a generalization of the $y^2 = g(x) h(x)$ example above,
where all possible factorizations are considered simultaneously.
The paper~\cite{BruinStoll2Desc} describes the procedure in detail.

\begin{Example} (See~\cite{BruinStollExp,BruinStoll2Desc})
  Among the small genus 2 curves, there are only $1492$ curves $C$
  without rational points and such that $\Sel^{(2)}(C) \neq \emptyset$.
  So 2-descent is a rather efficient tool in this case. Figure~1
  at the end of~\cite{BruinStoll2Desc} shows that this is still mostly
  true also for larger coefficients.
\end{Example}


\subsection{A Conjecture}

In the example above, we have seen that 2-descent shows that most of the
small curves without rational points really do not have rational points.
This makes it plausible that perhaps we can deal with the remaining
curves by an $n$-descent with a suitable $n > 2$. Unfortunately, the direct
computation of the relevant Selmer sets is infeasible. Still, we can
formulate the following conjecture. In~\cite{StollFiniteDescent}, we argue
that there are good reasons for it to hold.

\begin{Conjecture} \label{Conj1}
  If $C(\Q) = \emptyset$, then $\Sel^{(n)}(C) = \emptyset$ for some $n \ge 1$.
\end{Conjecture}

The case $n = 1$ is equivalent to checking for points everywhere locally on~$C$,
since $\operatorname{id}_C : C \to C$ is the only $1$-covering of~$C$.

\begin{Remarks} \strut
  \begin{enumerate}\addtolength{\itemsep}{2mm}
    \item In principle, $\Sel^{(n)}(C)$ is \emph{computable} for every~$n$.
          The conjecture therefore implies that ``$C(\Q) = \emptyset$?'' 
          is \emph{decidable}.
          (Search for points by day, compute $\Sel^{(n)}(C)$ by night.)
    \item The conjecture implies that the \emph{Brauer-Manin obstruction}
          is the \emph{only} obstruction against rational points on curves.
          (In fact, the conjecture is equivalent to this statement.)
          See~\cite{StollFiniteDescent} for details.
  \end{enumerate}
\end{Remarks}


\section{The Mordell-Weil Sieve}


\subsection{The Idea}

We now assume that we know explicit generators of the Mordell-Weil group~$J(\Q)$,
where $J$ is, as before, the Jacobian variety of the curve~$C$. By~\cite{Weil},
$J(\Q)$ is a finitely generated abelian group. It is clear that in the
diagram~\eqref{CovDiagram} we only need to consider those $n$-coverings~$X$
of~$J$ that actually have rational points. These $n$-coverings are of the form
\[ J \ni P \longmapsto nP + Q \in J \qquad \text{with $Q \in J(\Q)$;} \]
the shift $Q$ is only determined modulo~$nJ(\Q)$.

The set we are interested in is therefore
\[ \bigl\{Q + nJ(\Q) : \bigl(Q + nJ(\Q)\bigr) \cap \iota(C) \neq \emptyset\bigr\}
   \subset J(\Q)/nJ(\Q) \,.
\]
By the above, it contains the subset of the $n$-Selmer set of~$C$ that consists
of $n$-coverings of~$C$ with rational points.
We approximate the condition by testing it modulo~$p$ for a set of primes~$p$.

Let $S$ be a finite set of primes of good reduction for~$C$.
Consider the following diagram.
\begin{equation} \label{MWSDiagram}
   \xymatrix{ C(\Q) \ar[d] \ar[r]^{\iota} 
                & J(\Q) \ar[d] \ar[r] 
                & J(\Q)/nJ(\Q) \ar[d]^-{\beta} \\
              \displaystyle\prod_{p \in S} C(\F_p) \ar@<5pt>[r]^{\iota}
                 \ar@/_20pt/_{\alpha}[rr]
                 & \displaystyle\prod_{p \in S} J(\F_p) \ar@<5pt>[r]
                 & \displaystyle\prod_{p \in S} J(\F_p)/nJ(\F_p)
            }
\end{equation}
We \emph{can compute} the maps $\alpha$ and~$\beta$, since they only involve
finite objects. If their images do not intersect, then it follows that
$C(\Q) = \emptyset$. This method is known as the `Mordell-Weil Sieve'.
It was first used by Scharaschkin in his thesis~\cite{Scharaschkin}.
Flynn~\cite{Flynn} used it on a number of genus~2 curves.
In~\cite{BruinStollExp}, it is applied to the remaining undecided small
genus~2 curves.

\begin{Example}
  If $C$ is a small genus~2 curve without rational points, then either
  $C$ fails to have points everywhere locally, or $\Sel^{(2)}(C) = \emptyset$,
  or the Birch and Swinnerton-Dyer Conjecture for~$J$ implies that $C$
  has no embedding into~$J$ (this is needed for 42~curves),
  or else the Mordell-Weil Sieve with suitable parameters $S$ and~$n$
  proves that $C(\Q)$ is empty.
\end{Example}

In order to obtain this result, one needs a carefully optimized implementation
of the Mordell-Weil sieve. See~\cite{BruinStollMWS} for details. The
parameter the complexity depends on most sensitively is the rank~$r$
of~$J(\Q)$. If $r \le 3$, our implementation works quite well; it should
be mentioned that it uses not only information mod~$p$ for good primes~$p$,
but also information modulo powers of (small) primes, even when they are
primes of bad reduction. There are not yet enough worked examples where
the rank is larger than~3, so it is hard to say anything precise about
the performance of the algorithm in this case. At least there are isolated
examples that show that it can still work when $r$ is as large as~6.

\medskip

Poonen~\cite{PoonenHeur} shows that under reasonable assumptions, the
following should be true.

\begin{Conjecture}[Poonen Heuristic] \label{Conj2}
  If $C(\Q) = \emptyset$, then the maps $\alpha$ and~$\beta$ in
  Diagram~\eqref{MWSDiagram} above will have disjoint images when $n$
  and the set~$S$ are sufficiently large.
\end{Conjecture}

Conjecture~\ref{Conj2} implies Conjecture~\ref{Conj1} if we assume that
$\Sha(J/\Q)$ has no nontrivial infinitely divisible elements,
see~\cite{StollFiniteDescent}.


\subsection{Satisfying the Assumption on $J(\Q)$}

We are assuming here that we know explicit generators of~$J(\Q)$.
For the Mordell-Weil sieve as described above, if we use it to show that
$C$ has no rational points, it is actually sufficient to know generators
of a subgroup of finite index, if we can also show that the index is
coprime to~$n$. The latter is usually not so hard; see~\cite{FlynnSmart}.
To achieve the former, we can use $n$-descent again, but this time
on the Jacobian~$J$. This is feasible for hyperelliptic curves when $n = 2$
and in a few other rather special cases, 
see~\cite{Schaefer,PoonenSchaefer,Stoll2Desc,PSS}.
As with elliptic curves, large generators can be a problem, however.

\begin{Example} (See~\cite{BruinStollExp}.)
  For the small genus~2 curve
  \[ C : y^2 = -3\,x^6 + x^5 - 2\,x^4 - 2\,x^2 + 2\,x + 3 \,, \]
  the Mordell-Weil group $J(\Q)$ is infinite cyclic, generated by $[P_1 + P_2 - W]$,
  where the $x$-coordinates of $P_1$ and~$P_2$ are the roots of
  \[ x^2 + {\textstyle\frac{37482925498065820078878366248457300623}%
                {\mathstrut 34011049811816647384141492487717524243}}\,x 
         + {\textstyle\frac{581452628280824306698926561618393967033}%
                {\mathstrut 544176796989066358146263879803480387888}} \,,
  \]
  and $W$ is a canonical divisor.
\end{Example}

The bound on~$r$ obtained from 2-descent on~$J$ need not be tight.
The difference to the actual value of~$r$ comes from 2-torsion elements
of the Shafarevich-Tate group~$\Sha(J/\Q)$. In some cases, it is possible
to show that there are non-trivial such elements, thereby improving the
upper bound on the rank. Two techniques that have been suggested and also
used are \emph{visualization}~\cite{BruinFlynnVis}
and the Brauer-Manin obstruction on certain related
varieties~\cite{BBFL,FlynndP,LoganvanLuijk}.

For `reasonable' curves of genus~2, generators of a finite-index subgroup
of~$J(\Q)$ can usually be determined. For hyperelliptic curves of genus
at least~3, it may still be possible in many cases, but the situation is
already less favorable. Beyond hyperelliptic curves and variations on that
theme, descent calculations appear to be rather hopeless with the currently
available technology. There are some first attempts at 2-descent on
Jacobians of non-hyperelliptic curves of genus~3, however, so maybe the
situation will change at some point in the not-too-distant future.


\subsection{An Extension}

If we take $n$ in Diagram~\eqref{MWSDiagram} to be a multiple of a fixed
number~$N$, then we can restrict to a given \emph{coset~$X$ of~$N J(\Q)$}
(since this coset will be a union of cosets of~$n J(\Q)$). Therefore
the Mordell-Weil sieve computation gives us a way of proving that the
coset~$X$ does not meet~$\iota(C)$. In this case, there are no rational
points on~$C$ that are mapped into~$X$ under~$\iota$.

Conjecture~\ref{Conj1} can be extended to this situation.

\begin{Conjecture} \label{Conj3}
  Let $Q \in J(\Q)$.
  If $\bigl(Q + NJ(\Q)\bigr) \cap \iota(C) = \emptyset$,
  then there are $n \in N\Z$ and $S$ such that
  the Mordell-Weil sieve with these parameters proves this fact.
\end{Conjecture}

So if we can find an $N$ that \emph{separates} the rational points on~$C$,
i.e., such that the composition
$C(\Q) \stackrel{\iota}{\to} J(\Q) \to J(\Q)/NJ(\Q)$
is injective,
then we \emph{can effectively determine $C(\Q)$} if Conjecture~\ref{Conj3}
holds for~$C$. The procedure simply considers each coset of~$N J(\Q)$ in turn.
On the one hand, we run a search on~$C$ to find a rational point that maps
into the coset under consideration; on the other hand, we run the Mordell-Weil
sieve with the aim of proving that no such point exists. If Conjecture~\ref{Conj3}
holds, then one of these two computations has to produce a result. (In practice,
we just run the Mordell-Weil sieve. As long as the intersection of the images
of $\alpha$ and~$\beta$ is nonempty, we check the smallest representatives
in~$J(\Q)$ of the elements of the intersection whether they come from the curve.)


\section{Chabauty's Method}


\subsection{The Idea}

Chabauty~\cite{Chabauty} used this method to prove Mordell's Conjecture in
the case that the rank~$r$ of~$J(\Q)$ is smaller than the genus~$g$ of the
curve. The idea is to consider the $p$-adic points $J(\Q_p)$ as a $p$-adic
Lie group. The topological closure of~$J(\Q)$ then is a Lie subgroup of
dimension at most~$r$. One then expects that this subgroup of positive codimension
has only finitely many points of intersection with the analytic
curve~$\iota(C(\Q_p))$. This is what Chabauty proves. Later, the method
was taken up by Coleman~\cite{Coleman} who used it to deduce upper bounds
on the number of rational points on the curve. The method can also be
used to determine the set of rational points in certain cases,
see~\cite{FlynnChab,FPS,Wetherell} for early examples of this.
The book~\cite{CasselsFlynn} contains a description of the method when
$C$ has genus~2.

We now describe the setting more concretely. Let $p$ be a prime of good
reduction for~$C$ (this assumption simplifies things, but is not strictly
necessary). We denote by $\Omega_J^1(\Q_p)$ the $g$-dimensional $\Q_p$-vector
space of regular 1-forms on~$J$, and similarly for~$C$. Then $\iota$ induces
an isomorphism of $\Omega_J^1(\Q_p)$ and $\Omega_C^1(\Q_p)$ that is in fact
independent of our choice of the embedding~$\iota$.

The $p$-adic logarithm on~$J$ is a continuous group homomorphism
\[ \log : J(\Q_p) \To T_0 J(\Q_p) = \Omega_J^1(\Q_p)^* \]
whose kernel consists of the elements of finite order. It induces a pairing
\[ \Omega_J^1(\Q_p) \times J(\Q_p) \To \Q_p\,, \qquad
   (\omega, R) \longmapsto \int_0^R \omega = \langle \omega, \log R \rangle
\]
that becomes perfect if we replace $J(\Q_p)$ by $J(\Q_p)^0 \otimes_{\Z_p} \Q_p$
(where $J(\Q_p)^0$ is a sufficiently small neighborhood of the identity).
Since
\[ \rank \,J(\Q) = r < g = \dim_{\Q_p} \Omega_J^1(\Q_p) \,, \]
there is a differential
\[ 0 \neq \omega_p \in \Omega_C^1(\Q_p) \cong \Omega_J^1(\Q_p) \]
that kills $J(\Q) \subset J(\Q_p)$.

Let $P_0 \in C(\Q)$ be used as the base-point for the embedding~$\iota$.
Then the above implies that every point $P \in C(\Q)$ must satisfy
\[ \lambda(P) = \int_{P_0}^P \omega_p = 0 \,. \]
The function~$\lambda$ is a $p$-adic analytic function on~$C(\Q_p)$. On each
residue class mod~$p$, it can be represented by an explicit converging
power series. This makes it possible to bound the number of zeros of~$\lambda$
on such a residue class. If we find the same number of rational points within
the residue class, then we know that we have found them all. Some of the
zeros of~$\lambda$ may occur at transcendental points, however; in this case
the upper bound on the number of points is not tight. We can use information
from the Mordell-Weil sieve to rule out the spurious zeros; see~\cite{PSS}
for some examples.


\subsection{Combination with the Mordell-Weil Sieve}

We can also switch the roles of Chabauty's method and the Mordell-Weil sieve
and use Chabauty's method in a helping function. We still need to assume
that $r < g$. The idea is to use Chabauty's approach to find a
\emph{separating} number~$N$. The key to this is the following result
(see for example~\cite{StollChab}).

\begin{Theorem}
  Suppose that $p$ is a prime of good reduction for~$C$ and that
  $\omega_p \in \Omega_C^1(\Q_p)$ is a differential that kills the
  Mordell-Weil group~$J(\Q)$. We can assume that $\omega_p$ is scaled
  so that it has a well-defined reduction $\bar{\omega}_p \neq 0$ mod~$p$.
  If $\bar{\omega}_p$ does not vanish on~$C(\F_p)$ and $p > 2$,
  then each residue class mod~$p$ on~$C$ contains at most one rational point.
\end{Theorem}

In this case, the number $N = \#J(\F_p)$ is separating, since we know that
the map $C(\Q) \to C(\F_p)$ is injective (this is the statement of the theorem),
that $\iota : C(\F_p) \to J(\F_p)$ is injective, and that $J(\Q)/NJ(\Q)$
maps to~$J(\F_p)$.

Heuristic considerations indicate that the theorem applies for a set of primes~$p$
of positive density whenever $r < g$ and $J$ is simple. (If $J$ splits, we
can use one of the factors of~$J$ to do a similar computation.)

The most accessible case is when $g = 2$, since then we have a good chance
to determine~$J(\Q)$. The `Chabauty condition' $r < g$ then reduces to $r = 1$.
(When $r = 0$, the group $J(\Q)$ is finite, and we can easily find its
intersection with~$\iota(C)$, so this case is essentially trivial.)
In this case, the differentials $\bar{\omega}_p$ can be computed very
easily, and we quickly find suitable primes~$p$. The search for a suitable
separating number~$N$ can be integrated with the Mordell-Weil sieve computation.
This leads to a very efficient implementation that determines $C(\Q)$ quite fast
for genus~2 curves~$C$ such that $\rank\, J(\Q) = 1$.
See~\cite{BruinStollMWS} for a discussion.

\begin{Example} (See~\cite{StollDensity})
  For the $46\,436$ small genus 2 curves with rational points and such that
  $r = 1$, we determined $C(\Q)$. This computation takes about 8--9 hours
  on current hardware (as of~2009).
\end{Example}


\section{Some Odds and Ends}

In this section, we collect some remarks on extensions and variants of the
methods discussed above, and on some other approaches.

\subsection{Larger Rank}

When $r \ge g$, we can still use the Mordell-Weil Sieve
to show that we know all rational points up to very large height.
For this, we increase $n$ in the sieving procedure until we can prove
that none of the remaining cosets of~$n J(\Q)$ contains a point on~$C$
of height smaller than a given bound, except for the points we know.
Once the bulk of the computation is done, we can increase the height
bound without much extra cost.

If the desired height bound is not so large, it may be more efficient to
use lattice point enumeration. For this, we use the fact that the torsion-free
quotient $J(\Q)/J(\Q)_{\tors} \cong \Z^r$ is endowed with a positive definite
quadratic form~$\hat{h}$, the \emph{canonical height}. The height bound
for $P \in C(\Q)$ translates into a bound for $\hat{h}(\iota(P))$; we can
then enumerate all points in $J(\Q)$ (mod torsion) up to that height bound
and check if they come from the curve.

\begin{Example} (See~\cite{StollDensity})
  Unless there are points of height $> 10^{100}$,
  the largest point on a small genus~2 curve has height (i.e., maximum
  absolute value of the numerator and denominator of its $x$-coordinate)
  $209\,040$.
\end{Example}

For these applications, it is not enough to know generators of a finite-index
subgroup of~$J(\Q)$. We really need to know generators of the full
Mordell-Weil group. The `saturation' step from a finite-index subgroup
to the full group requires the computation of canonical heights on~$J$,
and we need a bound for the difference between the canonical height and
a suitable `naive height'. So far, the necessary theory and algorithms
only exist for $g \le 2$~\cite{FlynnHt,FlynnSmart,StollHt1,StollHt2}.
Therefore, we are currently limited to curves
of genus~2. There is current work that aims at extending the tools so that
they can also be used for higher-genus hyperelliptic curves, so we may
soon be able to deal with a larger class of curves.


\subsection{Integral Points}

If we can determine the set of rational points on~$C$, we obviously have
also found the \emph{integral} points. However, we can determine the set
of integral points in some cases even when we are not able to find~$C(\Q)$.
For example, if $C$ is hyperelliptic, we can compute bounds for
integral points using \emph{Baker's method} of `Linear forms in logarithms'.
The currently best results in this direction~\cite{BMSST} lead to bound
of a flavor like $|x| < 10^{10^{600}}$.

If we know \emph{generators} of $J(\Q)$, we can use the Mordell-Weil sieve
as explained in the previous subsection
to prove that there are no unknown rational points below that bound.
(The bound for $\hat{h}$ is something like $10^k$ with $k$ of the order
of several hundred or a couple of thousand. This is within reach of our
current implementation of the Mordell-Weil sieve method. See~\cite{BMSST}
or~\cite{BruinStollMWS} for details.)
It follows that we already know all the integral points on~$C$.

\begin{Example} (See~\cite{BMSST})
  The integral solutions to
  \[ \binom{y}{2} = \binom{x}{5} \]
  have $x \in \{0,1,2,3,4,5,6,7,15,19\}$.
\end{Example}

Since we need to know generators of the full Mordell-Weil group for this
application, the remarks made at the end of the previous subsection also
apply here. In particular, we are currently restricted to curves of
genus~2.


\subsection{Genus Larger Than Two}

We have seen that there are methods available that allow us to find out
a lot about the rational points on a given curve \emph{of genus~2}.
When the genus is larger, a number of difficulties arise.

If $C$ is hyperelliptic (or perhaps of some other rather special form),
it is still possible to do 2-descent on~$C$ and (to a certain extent) on~$J$.
For other curves, there is so far no feasible way to obtain provable
upper bounds on the rank of~$J(\Q)$. If we are willing to assume the
Birch and Swinnerton-Dyer conjecture for~$J$ (together with some related
conjectures on L-series) and the conductor of~$J$ is not too large,
then we can use Tim Dokchitser's code~\cite{Dokchitser} to compute
(an upper bound for) the order of vanishing of $L(s,J)$ at~$s=1$, which
gives a conditional upper bound on~$r$. We may then be able to find
a set of generators of a finite-index subgroup of~$J(\Q)$; this
suffices to apply Chabauty's method or its combination with the
Mordell-Weil sieve.

Another difficulty is the missing explicit theory of heights. This
prevents us from obtaining generators of the full Mordell-Weil group
(or rather, it prevents us from showing that we actually have generators).
This means that we cannot use the techniques described earlier in this
section.

Here are some examples that show what can still be done with curves
of genus at least~$3$.

\begin{Example} (See~\cite{PSS})
  In the course of solving $x^2 + y^3 = z^7$ in coprime integers, one has
  to determine the set of rational points on certain twists of the Klein Quartic.
  These are rather special non-hyperelliptic curves of genus~$3$.
  2-Descent on~$J$ is possible here; Chabauty and Mordell-Weil sieve
  techniques are successful.
\end{Example}

\begin{Example} (See~\cite{StollDyn})
  The curve $X_0^{\text{dyn}}(6)$ classifying 6-cycles under $x \mapsto x^2 + c$
  has genus 4. Assuming the Birch and Swinnerton-Dyer conjecture for its Jacobian,
  we can show that $r = 3$. We can then apply Chabauty's method to determine 
  $X_0^{\text{dyn}}(6)(\Q)$. It follows that there are no 6-cycles consisting
  of rational numbers (under the assumptions made).
\end{Example}

\begin{Example} (See~\cite{SiksekStoll})
  What are the arithmetic progressions in coprime integers that have the
  form $(a^2, b^2, c^2, d^5)$? This question leads to a number of hyperelliptic
  curves of genus~4; every solution to the original question gives rise to
  a rational point on one of these curves. There are three essentially different
  curves. For two of them, the 2-Selmer set turns out to be empty. For the
  last one, a 2-descent on its Jacobian is possible and shows that the rank is~2.
  Chabauty, combined with a little Mordell-Weil sieve information, then
  succeeds in showing that there are no unexpected points. This finally
  proves that the only arithmetic progression of the desired form is the
  trivial one, $(1,1,1,1)$.
\end{Example}


\subsection{The Method of Dem'yanenko-Manin}

The method of Dem'yanenko-Manin is an alternative method that can be used
to determine~$C(\Q)$ in some cases. When it applies, it gives an effective
bound on the height of the rational points on~$C$ (and not just on their
number, as is the case with Chabauty's method).

The requirement here is that we have $m$ independent morphisms $C \to A$,
where $A$ is some abelian variety and $m > \rank\, A(\Q)$. The idea is
that the images of points on~$C$ under these independent morphisms want
to be independent in~$A(\Q)$, but there is not enough room for them to
be independent. This leads to a bound on the height of the points.

If one looks at more or less `random' curves~$C$ that have two independent
maps to an elliptic curve~$E$, say, then the two images on~$E$ of a rational
point on~$C$ usually \emph{are} independent in~$E(\Q)$, invalidating the assumption.
So the method appears to be of fairly limited applicability.

There are cases, however, when the method can be used with profit.
In~\cite{Demjanenko}, it is applied to certain twists of the Fermat quartic
that have two independent maps to an elliptic curve. However, as Serre comments
in~\cite{SerreMW}, it is hard to find nontrivial examples.
See~\cite{GirardKulesz} for a more recent variation on this theme.

In~\cite{Manin}, Manin makes use of the growing number
of degeneracy maps $X_0(p^n) \to X_0(p)$ in order to show that for any
given prime~$p$, the power of~$p$ that divides the order of a rational torsion
point on an elliptic curve over~$\Q$ (or over any fixed number field)
is bounded.


\subsection{Covering Collections and Elliptic Curve Chabauty}

The Descent Theorem~\ref{Thm:Descent} tells us that we obtain all rational
points on a given curve~$C$ from the rational points on the various twists~$D_\xi$
of a covering of~$C$. If we can find this collection of twists explicitly
(this is sometimes called a \emph{covering collection} for~$C$), then we
can attempt to determine their sets of rational points instead of directly
trying to find~$C(\Q)$. This can be helpful when the rank of~$J(\Q)$ is too
large to apply Chabauty's method on~$C$, since the ranks associated to the
curves~$D_\xi$ may well be sufficiently small.

The downside of this approach is that the covering curves have larger
genus than~$C$, and so the methods described here are usually not applicable.
In some cases, the curves~$D_\xi$ map to other curves of low genus. If we
can find their rational points, we can determine those on~$D_\xi$. A very
useful variant arises when the target is an elliptic curve~$E$; the map may be
defined over some number field~$K$. The images of rational points on~$D_\xi$
then satisfy some additional constraints. This can be used to find these
images by a variant of Chabauty's method (applied to the restriction of
scalars of~$E$ from $K$ down to~$\Q$) when the rank of~$E(K)$ is less than
the degree of~$K$. This is known as \emph{Elliptic curve Chabauty};
see~\cite{Bruin:Thesis,Bruin:Ellchab,BruinElkies,BruinFlynn,Wetherell}
for details and examples.


\section{Concluding Remarks}

The last ten or fifteen years have seen tremendous progress in our ability
to determine the set of rational points on curves of higher genus, in particular
on curves of genus~2. Given a curve~$C$ of genus~2 over~$\Q$, we can now do the
following.

\begin{itemize}\addtolength{\itemsep}{1mm}
  \item Search for rational points on~$C$.
  \item Check if $C$ has points everywhere locally.
  \item Perform a 2-descent on~$C$, thus possibly showing that $C(\Q)$ is empty.
  \item Perform a 2-descent on~$J$, the Jacobian of~$C$, thus obtaining
        an upper bound on $r = \rank\, J(\Q)$.
  \item Search for rational points on~$J$, thus obtaining a lower bound on~$r$.
  \item Find generators of a finite-index subgroup of~$J(\Q)$ if both bounds
        agree.
  \item Compute canonical heights on~$J$.
  \item Find generators of~$J(\Q)$ if generators of a finite-index subgroup
        are known, assuming that the bound for the difference between naive
        and canonical height is not too large.
  \item If $r \le 1$, determine $C(\Q)$ using a combination
        of the Mordell-Weil sieve and Chabauty's method. (Termination of
        this is conditional on Conjecture~\ref{Conj3}, but if the
        computation terminates, which is always the case in practice,
        the result is provably correct.)
  \item If $r \ge 2$ and $J(\Q)$ is known, find all rational points on~$C$ up to
        very large height.
  \item If $J(\Q)$ is known, find all integral points on~$C$.
\end{itemize}

From a practical point of view, what is missing to make this really satisfying
is a way of determining a separating~$N$ when $r \ge 2$ (without previous
knowledge what $C(\Q)$ is). If a separating~$N$ can be found, then the same
approach as used when combining the Mordell-Weil sieve with Chabauty's
method will enable us to determine~$C(\Q)$.

From a theoretical point of view, we would like to have a proof of 
Conjecture~\ref{Conj3}, since this will guarantee that our procedure terminates.
(For practical computations, we don't really care about a proof as long as
the computation terminates; the result will be correct in any case.)
The other theoretical gap is that it is still open whether the rank~$r$
can be found effectively. This is related to the finiteness of~$\Sha(J/\Q)$,
which is only known in very special cases.

For curves of higher genus than~2, some of the items on the list above
can still be done (in particular when $C$ is hyperelliptic), but we soon
reach a point where things become infeasible. However, I believe that this
is only a matter of complexity and not of principle: given sufficient
resources, we should be able to perform the same kind of computation also
with more general curves. (Of course, some theoretical work still has to
be done for this, like an extension of the explicit theory of heights that
we have at our disposal when the genus is~2.)

\medskip

Based on what we can actually do, on various heuristic considerations, and
on fairly extensive experimental data, I am convinced that it is actually
possible (in principle) to determine the set~$C(\Q)$ algorithmically,
when $C$ is a curve of genus~$\ge 2$. A complete proof of this statement
is likely to be quite far away still, but the progress that has been made
on the practical side in recent years is very encouraging.



\begin{thebibliography}{99}
  \frenchspacing

  \bibitem{BBFL}
    {\sc M.J. Bright, N. Bruin, E.V. Flynn, A. Logan:}
    {\it The Brauer-Manin obstruction and $\Sha[2]$,}
    LMS J. Comput. Math. {\bf 10} (2007), 354--377.

  \bibitem{Bruin:Thesis}
    {\sc N. Bruin:} 
    {\it Chabauty methods and covering techniques applied to generalized Fermat
    equations,} CWI Tract 133, 77 pages (2002).

  \bibitem{Bruin:Ellchab}
    {\sc N. Bruin:}
    {\it Chabauty methods using elliptic curves,}
    J. Reine Angew. Math. {\bf 562} (2003), 27--49.

  \bibitem{BruinElkies}
    {\sc N. Bruin, N.D. Elkies:}
    {\it Trinomials $ax\sp 7+bx+c$ and $ax\sp 8+bx+c$ with Galois groups of order
     168 and $8\cdot 168$,} in:
    {\it Algorithmic number theory, Sydney 2002,}
    Lecture Notes in Comput. Sci. {\bf 2369}, Springer, Berlin (2002),
    pp. 172--188.

  \bibitem{BruinFlynn}
    {\sc N. Bruin, E.V. Flynn:}
    {\it Towers of 2-covers of hyperelliptic curves,}
    Trans. Amer. Math. Soc. {\bf 357} (2005), 4329--4347.

  \bibitem{BruinFlynnVis}
    {\sc N. Bruin, E.V. Flynn:}
    {\it Exhibiting SHA$[2]$ on hyperelliptic Jacobians,}
    J. Number Theory {\bf 118} (2006), 266--291.

  \bibitem{BruinStollExp}
    {\sc N. Bruin, M. Stoll}:
    {\it Deciding existence of rational points on curves: an experiment,}
    Experiment. Math. {\bf 17} (2008), 181--189.

  \bibitem{BruinStoll2Desc}
    {\sc N. Bruin, M. Stoll:}
    {\it 2-cover descent on hyperelliptic curves},
    Math. Comp. {\bf 78} (2009), 2347--2370.

  \bibitem{BruinStollMWS}
    {\sc N. Bruin, M. Stoll:}
    {\it The Mordell-Weil sieve: Proving non-existence of rational points on 
    curves}, Preprint (2009), to appear in LMS J. Comput. Math. (2010).

  \bibitem{BMSST}
    {\sc Y. Bugeaud,  M. Mignotte, S. Siksek, M. Stoll, Sz. Tengely:}
    {\it Integral points on hyperelliptic curves,}
    Algebra Number Theory {\bf 2} (2008), 859--885.

  \bibitem{Cassels}
    {\sc J.W.S. Cassels:}
    {\it Second descents for elliptic curves,}
    J. reine angew. Math. {\bf 494} (1998), 101--127.

  \bibitem{CasselsFlynn}
    {\sc J.W.S. Cassels, E.V. Flynn:}
    {\it Prolegomena to a middlebrow arithmetic of curves of genus~2},
    London Math. Soc., Lecture Note Series {\bf 230},
    Cambridge Univ. Press, Cambridge, 1996.

  \bibitem{Chabauty}
    {\sc C. Chabauty:}
    {\it Sur les points rationnels des courbes alg\'ebriques de genre
    sup\'erieur \`a l'unit\'e,}
    C. R. Acad. Sci. Paris {\bf 212} (1941), 882--885.

  \bibitem{ChevalleyWeil}
    {\sc C. Chevalley, A. Weil:} {\it Un th\'eor\`eme d'arithm\'etique
    sur les courbes alg\'ebriques,} Comptes Rendus Hebdomadaires des
    S\'eances de l'Acad. des Sci., Paris {\bf 195} (1932), 570--572.

  \bibitem{Coleman}
    {\sc R.F. Coleman:} {\it Effective Chabauty,} Duke Math. J. {\bf 52} (1985),
    765--770.

  \bibitem{CFOSS}
    {\sc J.E. Cremona, T.A. Fisher, C. O'Neil, D. Simon, M. Stoll:}
    {\it Explicit $n$-descent on elliptic curves.}
    {\it I. Algebra,} J. reine angew. Math. {\bf 615} (2008), 121--155.
    {\it II. Geometry,} J. reine angew. Math. {\bf 632} (2009), 63--84.
    {\it III. Algorithms,} in preparation.

  \bibitem{CFS}
    {\sc J.E. Cremona, T.A. Fisher, M. Stoll:}
    {\it Minimisation and reduction of 2-, 3- and 4-coverings of elliptic curves,}
    Preprint (2009), arXiv:0908.1741v1 [math.NT].

  \bibitem{Demjanenko}
    {\sc V.A. Dem'janenko:}
    {\it Rational points of a class of algebraic curves} (Russian),
    Izv. Akad. Nauk SSSR Ser. Mat. {\bf 30} (1966), 1373--1396.

  \bibitem{Dokchitser}
    {\sc T. Dokchitser:} {\it Computing special values of motivic $L$-functions,}
    Experiment. Math. {\bf 13} (2004), 137--149.

  \bibitem{Faltings}
    {\sc G. Faltings:}
    {\it Endlichkeitss\"atze f\"ur abelsche Variet\"aten \"uber Zahlk\"orpern,}
    Invent. Math. {\bf 73} (1983), 349--366.

  \bibitem{FlynnHt}
    {\sc E.V. Flynn:} {\it An explicit theory of heights,}
    Trans. Amer. Math. Soc. {\bf 347} (1995), 3003--3015.

  \bibitem{FlynnChab}
    {\sc E.V. Flynn:}
    {\it A flexible method for applying Chabauty's theorem,}
    Compositio Math. {\bf 105} (1997), 79--94.

  \bibitem{Flynn}
    {\sc E.V. Flynn:}
    {\it The Hasse Principle and the Brauer-Manin obstruction for curves,}
    Manuscripta Math. {\bf 115} (2004), 437--466.

  \bibitem{FlynndP}
    {\sc E.V. Flynn:}
    {\it Homogeneous spaces and degree 4 del Pezzo surfaces,}
    Manuscripta Math. {\bf 129} (2009), 369–380.

  \bibitem{FPS}
    {\sc E.V. Flynn, B. Poonen, E.F. Schaefer:}
    {\it Cycles of quadratic polynomials and rational points on a genus-2 curve,}
    Duke Math. J. {\bf 90} (1997), 435--463.

  \bibitem{FlynnSmart}
    {\sc E.V. Flynn, N.P. Smart:}
    {\it Canonical heights on the Jacobians of curves of genus $2$
    and the infinite descent,}
    Acta Arith. {\bf 79} (1997), 333--352.

  \bibitem{GirardKulesz}
    {\sc M. Girard, L. Kulesz:}
    {\it Computation of sets of rational points of genus-3 curves via 
    the Dem'janenko-Manin method,}
    LMS J. Comput. Math. {\bf 8} (2005), 267--300.

  \bibitem{Ih}
    {\sc Su-Ion Ih:}
    {\it Height uniformity for algebraic points on curves,}
    Compositio Math. {\bf 134} (2002), 35--57.

  \bibitem{Kolyvagin}
    {\sc V.A. Kolyvagin:} {\it Finiteness of $E(\Q)$ and $\Sha(E,\Q)$ for a 
    subclass of Weil curves,} Izv. Akad. Nauk SSSR Ser. Mat., Vol.~{\bf 52} (1988),
    522--540.

  \bibitem{LLL}
    {\sc A.K. Lenstra, H.W. Lenstra, Jr., L. Lov\'asz:}
    {\it Factoring polynomials with rational coefficients,}
    Math. Ann. {\bf 261} (1982), 515--534.

  \bibitem{LoganvanLuijk}
    {\sc A. Logan, R. van Luijk:}
    {\it Nontrivial elements of Sha explained through K3 surfaces,}
    Math. Comp. {\bf 78} (2009), 441--483.

  \bibitem{Manin}
    {\sc Y. Manin:}
    {\it The $p$-torsion of elliptic curves is uniformly bounded} (Russian),
    Izv. Akad. Nauk SSSR Ser. Mat. {\bf 33} (1969), 459--465.

  \bibitem{MSS}
    {\sc J.R. Merriman, S. Siksek, N.P. Smart:}
    {\it Explicit $4$-descents on an elliptic curve,}
    Acta Arith. {\bf 77} (1996), 385--404.

  \bibitem{Mordell}
    {\sc L.J. Mordell:}
    {\it On the rational solutions of the indeterminate equations
     of the 3rd and 4th degrees,}
    Proc. Camb. Phil. Soc. {\bf 21} (1922), 179--192.

  \bibitem{PoonenHeur}
    {\sc B. Poonen:}
    {\it Heuristics for the Brauer-Manin obstruction for curves,}
    Experiment. Math. {\bf 15} (2006), 415--420.

  \bibitem{PoonenSchaefer}
    {\sc B. Poonen, E.F. Schaefer:} 
    {\it Explicit descent for Jacobians of cyclic covers of the projective line,} 
    J. reine angew. Math. {\bf 488} (1997), 141--188.

  \bibitem{PSS}
    {\sc B. Poonen, E.F. Schaefer, M. Stoll:}
    {\it Twists of $X(7)$ and primitive solutions to $x^2+y^3=z^7$,}
    Duke Math. J. {\bf 137} (2007), 103--158.

  \bibitem{PoonenStoll}
    {\sc B. Poonen, M. Stoll:} 
    {\it A local-global principle for densities,}
    in: {\sc Scott D. Ahlgren} (ed.) et al.: 
    {\it Topics in number theory. In honor of B. Gordon and S. Chowla.}
    Kluwer Academic Publishers, Dordrecht. 
    Math. Appl., Dordr. {\bf 467} (1999), 241--244.

  \bibitem{Schaefer}
    {\sc E.F. Schaefer:}
    {\it Computing a Selmer group of a Jacobian using functions on the curve,}
    Math. Ann. {\bf 310} (1998), 447--471.

  \bibitem{Scharaschkin}
    {\sc V. Scharaschkin:}
    {\it Local-global problems and the Brauer-Manin obstruction,}
    Ph.D. thesis, University of Michigan (1999).

  \bibitem{SerreGCFT}
    {\sc J.-P. Serre:} {\it Algebraic groups and class fields}, Springer
    GTM {\bf 117}, Springer Verlag, 1988.

  \bibitem{SerreMW}
    {\sc J.-P. Serre:}
    {\it Lectures on the Mordell-Weil theorem.}
    Translated from the French and edited by Martin Brown from notes by Michel 
    Waldschmidt. Aspects of Mathematics, E15. Friedr. Vieweg \& Sohn,
    Braunschweig, 1989.

  \bibitem{SiksekStoll}
    {\sc S. Siksek, M. Stoll:}
    {\it On a problem of Hajdu and Tengely,}
    Preprint (2009), arXiv:0912.2670v1 [math.NT].

  \bibitem{Simon}
    {\sc D. Simon:}
    {\it Solving quadratic equations using reduced unimodular quadratic forms,}
    Math. Comp. {\bf 74} (2005), 1531--1543.

  \bibitem{Stamminger}
    {\sc S. Stamminger:} {\it Explicit 8-descent on elliptic curves,}
    PhD thesis, International University Bremen (2005).

  \bibitem{StollHt1}
    {\sc M. Stoll:} {\it On the height constant for curves of genus two,}
    Acta Arith. {\bf 90} (1999), 183--201.

  \bibitem{Stoll2Desc}
    {\sc M. Stoll:} {\it Implementing 2-descent on Jacobians of hyperelliptic
    curves,} Acta Arith. {\bf 98} (2001), 245--277.

  \bibitem{StollHt2}
    {\sc M. Stoll:} {\it On the height constant for curves of genus two, II,}
    Acta Arith. {\bf 104} (2002), 165--182.

  \bibitem{StollIHP}
    {\sc M. Stoll:}
    {\it Descent on Elliptic Curves.}
    Short Course taught at IHP in Paris, October 2004.
    arXiv:math/0611694v1 [math.NT].

  \bibitem{StollChab}
    {\sc M. Stoll:}
    {\it Independence of rational points on twists of a given curve,}
    Compositio Math. {\bf 142} (2006), 1201--1214.

  \bibitem{StollFiniteDescent}
    {\sc M. Stoll:}
    {\it Finite descent obstructions and rational points on curves,}
    Algebra Number Theory {\bf 1} (2007), 349--391.

  \bibitem{StollDyn}
    {\sc M. Stoll:}
    {\it Rational 6-cycles under iteration of quadratic polynomials,}
    London Math. Soc. J. Comput. Math. {\bf 11} (2008), 367--380.

  \bibitem{StollDensity}
    {\sc M. Stoll:}
    {\it On the average number of rational points on curves of genus 2,}
    Preprint (2009), arXiv:0902.4165v1 [math.NT].

  \bibitem{Ratpoints}
    {\sc M. Stoll:}
    {\it Documentation for the ratpoints program,}
    Manuscript (2009), arXiv:0803.3165 [math.NT].

  \bibitem{Weil}
    {\sc A. Weil:} 
    {\it L'arithm\'etique sur les courbes alg\'ebriques,}
    Acta Math. {\bf 52} (1929), 281--315.

  \bibitem{Wetherell}
    {\sc J.L. Wetherell:} {\it Bounding the number of rational points on
    certain curves of high rank,} Ph.D. thesis, University of California
    (1997).

\end{thebibliography}
\end{document}